\documentclass{article}
\usepackage{amsthm, amsmath, amsfonts, amssymb}
\usepackage{charter, graphicx, mathrsfs}
\usepackage{authordate1-4}

\begin{document}

\renewcommand{\qedsymbol}{$\blacksquare$}
\newtheorem{thm}{Theorem} [section]     
\newtheorem{cor}[thm]{Corollary} 
\newtheorem{lem}[thm]{Lemma} 
\newtheorem{prop}[thm]{Proposition} 
\newtheorem{rem}[thm]{Remark}

\theoremstyle{definition}
\newtheorem{dfn}[thm]{Definition}
\newtheorem{ex}[thm]{Example}

\title{Homotopies in Classical and Paraconsistent Modal Logics}
\author{\textbf{Can Ba\c{s}kent}}
\date{}

\maketitle

\section{Introduction}

Bisimulations and van Benthem's celebrated theorem provide a direct insight how truth preserving operations in modal logics work. Apart from bisimulations, there are several other operations in basic modal logic that preserve the truth \cite{bla1}.

However, there is a problem. Given a modal model, and several bisimilar copies of it, there is no method to compare or measure the differences between bisimilar models apart from the basic model theoretical methods (i.e. they are submodels of each other, for instance). A rather negative slogan for this issue is the following: \emph{Modal language cannot distinguish bisimilar models}.

Since the modal language \emph{cannot count or measure},  several extensions of the language has been proposed to tackle this issue such as hybrid logics and majority logics \cite{bla,fin0}. However, such extensions are language based and introduce a non-natural, and sometimes counter-intuitive operators to the language, and often criticized as being \emph{ad-hoc}.

In this work, we will focus on the truth preserving operations rather than extending the modal language. In other words, we will ask the following question: \emph{Is there a truth preserving \textrm{natural} modal operation that can also distinguish the models it generates, and even compares them?} We argue that such an operation exist, and the answer to that question is positive.

However, to conceptualize our concerns and questions, we need to be careful at picking the correct modal logical framework. Kripkean models, in this respect, are criticized as they are overly simplistic and can overshadow some mathematical properties that can be apparent in some other modal models. Therefore, in this paper, we will concentrate on topological models for modal logics. On the other hand, note that topological models historically precede Kripke models, and are mathematically more complex allowing us to express variety of ideas within modal logic. Therefore, they can provide us with much stronger and richer structure of which we can take advantage. In this paper, we utilize a rather elementary concept from topology. We first introduce homemorphisms, and then homotopies to the modal logical framework, and show the immediate invariance results. On the other hand, from an application oriented point of view, we also have some applications to illustrate how our constructions can be useful.

In a previous work, homeomorphisms and homotopies were introduced to the context of nonclassical modal logics \cite{bas3}. In this work, our goal is to extend such results to both nonclassical and classical cases and see how homotopies can be defined in classical modal case. Before moving on, let us make it clear that what we mean by \emph{nonclassical} is either paraconsistent or paracomplete logical systems.

\subsection{What is a Topology?}

The history of the topological semantics of (modal) logics can be traced back to early 1920s making it the first semantics for variety of modal logics \cite{gold2}. The major revival of the topological semantics of modal logics and its connections with algebras, however, is due to McKinsey and Tarski \cite{mck,mck1}. In this section, we will briefly mention the basics of topological semantics in order to be able build our future constructions. We will give two equivalent definitions of topological spaces here for our purposes.

\begin{dfn}{\label{open-top}}
The structure $\langle S, \tau \rangle$ is called a topological space if it satisfies the following conditions.
\begin{enumerate}
\item $S \in \tau$ and $\emptyset \in \tau$
\item $\tau$ is closed under arbitrary unions and finite intersections
\end{enumerate}
\end{dfn}

\begin{dfn}{\label{close-top}}
The structure $\langle S, \sigma \rangle$ is called a topological space if it satisfies the following conditions.
\begin{enumerate}
\item $S \in \sigma$ and $\emptyset \in \sigma$
\item $\sigma$ is closed under finite unions and arbitrary intersections
\end{enumerate}
\end{dfn}

Collections $\tau$ and $\sigma$ are called topologies. The elements of $\tau$ are called \emph{open} sets whereas the elements of $\sigma$ are called \emph{closed} sets. Therefore, a set is open if its complement in the same topology is a closed set and vice versa. A function is called \emph{continuous} if the inverse image of an open (respectively, closed) set is open (respectively, closed), and a function is called \emph{open} if the image of an open (respectively, closed) set is open (respectively, closed). Moreover, two topological spaces are called \emph{homeomorphic} if there is a continuous bijection from one to the other with a continuous inverse. Moreover, two continuous functions are called \emph{homotopic} if there is a continuous deformation between the two. Homotopy is then an equivalence relation and, it gives rise to homotopy groups which is a foundational subject in algebraic topology.

\subsection{What is Paraconsistency?}

An easy and immediate semantics of paraconsistent/paracomplete logics can be given by using topologies. For this reason, it is helpful to remember some basics of paraconsistency.

First, note that deductive explosion describes the situation where any formula can be deduced from an inconsistent set of formulae, i.e. for all formulae $\varphi$ and $\psi$, we have $\{ \varphi, \neg \varphi \} \vdash \psi$, where $\vdash$ denotes logical consequence relation. In this respect, both ``classical'' and intuitionistic logics are known to be explosive. Paraconsistent logic, on the other hand, is the umbrella term for logical systems where the logical consequence relation $\vdash$ is \emph{not}  explosive \cite{pri0}. Variety of philosophical and logical objections can be raised against paraconsistency, and almost all of these objections can be defended in a rigorous fashion. We will not here be concerned about the philosophical implications of it, yet we refer the reader to the following for a comprehensive defense of paraconsistency with a variety of well-structured applications \cite{pri1}.

Use of topological semantics for paraconsistent logics is not new. To our knowledge, the earliest work discussing the connection between inconsistency and topology goes back to Goodman \cite{good}\footnote{Thanks to Chris Mortensen for this remark. Even if Goodman's paper appeared in 1981, the work had been carried out around 1978. In his paper, Goodman indicted that the results were based on an early that appeared in 1978 only as an abstract.}. In his paper, Goodman discussed ``pseudo-complements'' in a lattice theoretical setting and called the topological system he obtains ``anti-intuitionistic logic''. In a recent work, Priest discussed the dual of the intuitionistic negation operator and considered that operator in topological framework \cite{pri2}. Similarly, Mortensen discussed topological separation principles from a paraconsistent and paracomplete point of view and investigated the theories in such spaces \cite{mor}. Similar approaches from modal perspective was discussed by B\'{e}ziau, too \cite{bez}.

\subsection{Semantics}

In our setting, we denote set of propositional variables with $P$. We use the language of propositional modal logic with the modality $\Box$, and we will define the dual $\Diamond$ in the usual sense, and construct the language of the basic unimodal logic recursively in the standard fashion.

In topological semantics, the modal operator for necessitation corresponds to the topological \emph{interior} operator $\mathsf{Int}$ where $\mathsf{Int}(O)$ is the largest open set contained in $O$. Furthermore, one can dually associate the topological closure operator $\mathsf{Clo}$ with the possibility modal operator $\Diamond$ where the closure $\mathsf{Clo}(O)$ of a given set $O$ is the smallest closed set that contains $O$.

Before connecting topology and modal logic, let us set a piece of notation and terminology. The extension, i.e. the points at which the formula is satisfied, of a formula $\varphi$ in the model $M$ will be denoted as $[\varphi]^{M}$. We omit the superscript if the model we are working with is obvious. Moreover, by a \emph{theory}, we mean a deductively closed set of formulae.

The extensions of Boolean cases are obvious. However, the extension of a modal formula $\Box \varphi$ is then associated with an open set in the topological system. Thus, we have $[\Box \varphi] = \mathsf{Int}([\varphi])$. Similarly, we put $[\Diamond \varphi] = \mathsf{Clo}([\varphi])$. This means that in the basic setting, topological entities such as open or closed sets appear only with modalities. 

However, we can take one step further and suggest that extension of \emph{any} propositional variable be an open set \cite{mor,min}. In that setting, conjunction and disjunction works fine for finite intersections and unions. Nevertheless, the negation can be difficult as the complement of an open set is not generally an open set, thus may not be the extension of a formula in the language. For this reason, we need to use a new negation symbol $\dot{\sim}$ that returns the open complement (interior of the complement) of a given set. A similar idea can also be applied to closed sets where we assume that the extension of any propositional variable will be a closed set. In order to be able to avoid a similar problem with the negation, we stipulate yet another negation operator which returns the closed complement (closure of the complement) of a given set. In this setting, we use the symbol $\sim$ that returns the closed complement of a given set. Under these assumptions and recalling Definitions~\ref{open-top}~\&~\ref{close-top}, it is easy to observe the following \cite{mor}.
\begin{itemize}
\item In the topology of open sets $\tau$, any theory that includes the theory of the propositions that are true at the boundary is incomplete.
\item In the topology of closed sets $\sigma$, any theory that includes the boundary points will be inconsistent. 
\end{itemize}

An immediate observation yields that since extensions of all formulae in $\sigma$ (respectively in $\tau$) are closed (respectively, open), the topologies which are obtained in both paraconsistent (and paracomplete) logics are discrete. We already made the following simple connection \cite{bas3}. For a given model $M$, let $|M|$ denote the size of $M$'s carrier set.

\begin{prop}
Let $M_1$ and $M_2$ be paraconsistent and paracomplete topological models respectively. If $|M_1| = |M_2|$, then there is a homeomorphism from a paraconsistent topological model to the paracomplete one, and vice versa.
\end{prop}

\section{Homotopies}

Let us clarify an important point. We stipulate that the extension of any formula is closed to obtain a paraconsistent system, and stipulate that to be open to obtain a paracomplete set. In other words, we do not \emph{mix} such systems. On the other hand, note that in the classical case, open or closed sets appear only under the presence of modal operators, and they may appear together.

Now, in this section, we will start off with the easier case and consider paraconsistent topological spaces. Second, we will extend our results to classical case.

\subsection{Paraconsistent Case}

A recent research program that considers topological modal logics with continuous functions were discussed in an early work \cite{art1,kre}. An immediate theorem, which was stated and proved in variety of different work, would also work for paraconsistent logics \cite{kre}. Now, let us take two closed set topologies $\sigma$ and $\sigma'$ on a given set $S$ and a homeomorphism $f:\langle S, \sigma \rangle \rightarrow \langle S, \sigma' \rangle$. Akin to a previous theorem of Kremer and Mints, we have a simple way to associate the respective valuations between two models $M$ and $M'$ which respectively depend on $\sigma$ and $\sigma'$ so that we can have a truth preservation result. Therefore, define $V'(p)= f(V(p))$. Then, we have $M \models \varphi$ iff $M' \models \varphi$.

\begin{thm}{\label{cont-thm}}
Let $M = \langle S, \sigma, V \rangle$ and $M' = \langle S, \sigma', V' \rangle$ be two paraconsistent topological models with a homeomorphism $f$ from $\langle S, \sigma \rangle$ to $\langle S, \sigma' \rangle$. Define $V'(p)= f(V(p))$. Then $M \models \varphi$ iff $M' \models \varphi$ for all $\varphi$.
\end{thm}

\begin{proof}
The proof is by induction on the complexity of the formulae. See \cite{bas3} for the proof in non-classical case.
\end{proof}

Assuming that $f$ is a homeomorphism may seem a bit strong. We can then separate it into two chunks. One direction of the biconditional can be satisfied by continuity whereas the other direction is satisfied by the openness of $f$. 

\begin{cor}{\label{cor-open}}
 Let $M = \langle S, \sigma, V \rangle$ and $M' = \langle S, \sigma', V' \rangle$ be two paraconsistent topological models with a continuous $f$ from $\langle S, \sigma \rangle$ to $\langle S, \sigma' \rangle$. Define $V'(p)= f(V(p))$. Then $M \models \varphi$ implies $M' \models \varphi$ for all $\varphi$.
\end{cor}

\begin{cor}
 Let $M = \langle S, \sigma, V \rangle$ and $M' = \langle S, \sigma', V' \rangle$ be two paraconsistent topological models with an open $f$ from $\langle S, \sigma \rangle$ to $\langle S, \sigma' \rangle$. Define $V'(p)= f(V(p))$. Then $M' \models \varphi$ implies $M \models \varphi$ for all $\varphi$.
\end{cor}

Proofs of both corollaries depend on the fact that $\mathsf{Clo}$ operator commutes with continuous functions in one direction, and it commutes with open functions in the other direction. Furthermore, similar corollaries can be given for paracomplete frameworks as the $\mathsf{Int}$ operator also commutes in one direction under similar assumptions, and we leave it to the reader as well.

Furthermore, any topological operator that commutes with continuous, open and homeomorphic functions will reflect the same idea and preserve the truth\footnote{Thanks to Chris Mortensen for pointing this out.}. Therefore, these results can easily be generalized.

We can now take one step further to discuss homotopies in paraconsistent topological modal models. To best of our knowledge, homotopies first introduced to (non-classical) modal logic in \cite{bas3}. Now, recall that a \emph{homotopy} is a description of how two continuous function from a topological space to another can be deformed to each other. We can now state the formal definition.

\begin{dfn} Let $S$ and $S'$ be two topological spaces with continuous functions $f, f': S \rightarrow S'$. A homotopy between $f$ and $f'$ is a continuous function $H : S \times [0, 1] \rightarrow S'$ such that if $s \in S$, then $H(s, 0) = f(s)$ and $H(s, 1) = g(s)$.
\end{dfn}

In other words, a homotopy between $f$ and $f'$ is a family of continuous functions $H_t: S \rightarrow S'$ such that for $t \in [0, 1]$ we have $H_0 = f$ and $H_1 = g$ and the map $t \mapsto H_t$ is continuous from $[0, 1]$ to the space of all continuous functions from $S$ to $S'$. Notice that homotopy relation is an equivalence relation. Thus, if $f$ and $f'$ are homotopic, we denote it with $f \approx f'$. We will now use homotopies to obtain a generalization of Theorem~\ref{cont-thm}.

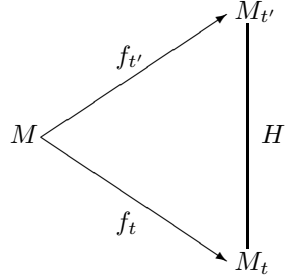
\begin{figure}[htp]
\begin{picture}(200,120)
\put(40,60){\mbox{$M$}}
\put(125,12){\mbox{$M_{t}$}}
\put(125,107){\mbox{$M_{t'}$}}
\put(80,27){\mbox{$f_{t}$}}
\put(80,90){\mbox{$f_{t'}$}}
\put(52,62){\vector(3,2){70}}
\put(52,62){\vector(3,-2){70}}
\put(130,105){\line(0,-1){85}}
\put(135,60){\mbox{$H$}}
\end{picture}
\caption{Homotopic Models}
\end{figure}

\begin{dfn}
Given a model $M = \langle S, \sigma, V \rangle$, we call the family of models $\{ M_t = \langle S_{t} \subseteq S, \sigma_{t}, V_t \rangle \}_{t \in [0, 1]}$ generated by $M$ and homotopic functions  homotopic models. In the generation, we put $V_t = f_t(V)$. 
\end{dfn}

\begin{thm}{\label{homotopy}}
Given two topological paraconsistent models $M = \langle S, \sigma, V \rangle$ and $M' = \langle S', \sigma', V' \rangle$ with two continuous functions $f, f' : S \rightarrow S'$ both of which respect the valuation: $V' = f(V)=f'(V)$. If there is a homotopy $H$ between $f$ and $f'$, then homotopic models satisfy the same formulae.
\end{thm}

\begin{proof}
To make the proof a bit more readable, note that the object whose names have a prime$'$, are in the range of the functions.

To make the proof go easily, we will assume that continuous functions are onto. If not, we can easily rearrange the range in such a way that it will be.

Let $M_{t} = \langle S'_{t}, \sigma_{t}, V_{t} \rangle$ and $M_{t'} = \langle S'_{t'}, \sigma_{t'}, V_{t'} \rangle$ be homotopic models with continuous $f_{t}, f_{t'}: S \rightarrow S'$. Observe that $M_{t} = \langle S'_{t}, \sigma_{t}, f_{t}(V) \rangle$, or equally $M_{t} = \langle S'_{t}, \sigma_{t}, H_{t}(V) \rangle$ for homotopy $H$.

Take a point $s_{t}$ such that $f_{t}(s)$ for some $s \in S$. 

Then, take $\langle S'_{t}, \sigma_{t}, H_{t}(V) \rangle, f_{t}(s) \models \varphi$ for arbitrary $\varphi$.   Since, $H$ is continuous on $t$ by some $h$, and by Corollary~\ref{cor-open}, we observe  $M_{t'} = \langle S_{t'}, \sigma_{t'}, H_{t'}(V) \rangle, s_{t'} \models \varphi$. In this case $s_{t'}$ exists as $H_{t}$ is continuos on $t$ and $s_{t'} = h(f_{t}(s)) = h(s_{t})$. Therefore, $M_{t} \models \varphi$ implies $M_{t'} \models \varphi$.

Notice that in this case, we did not need to present a proof on the complexity of $\varphi$. The reason for that is the fact that the extension of each formula is a closed set (since we are in a paraconsistent setting). 
\end{proof}

\begin{cor}
$M \models \varphi$ implies $M_{t} \models \varphi$, but not the other way around.
\end{cor}

\begin{cor}
In Theorem~\ref{homotopy}, if we take the cases for $t=0$ and $t=1$, we obtain Corollary~\ref{cor-open}.
\end{cor}

Notice that we have discussed the truth in the image sets that are obtained under $f, f', f_{t}, \dots$. Nevertheless, the converse can also be true, once the continuous functions have continuous inverses: this is exactly what is guaranteed by homeomorphisms. The corresponding notion at the level of homotopies is an \emph{isotopy}. An isotopy is a continuous transformation between homeomorphic functions. Thus, we have the following.

\begin{thm}
Given two topological paraconsistent models $M = \langle S, \sigma, V \rangle$ and $M' = \langle S', \sigma', V' \rangle$ with two homeomorphism $f, f' : S \rightarrow S'$ both of which respect the valuation: $V' = f(V)=f'(V)$. If there is an isotopy $H$ between $f$ and $f'$, then, for all $\varphi$, we have $$M_{t} \models \varphi \text{ iff } M \models \varphi \text{ iff } M' \models \varphi$$
\end{thm}

\begin{proof} Immediate, thus left to the reader. \end{proof}

What makes the non-classical case easy is the fact that the extension of each formula is an open or a closed set. Furthermore, a similar theorem can be stated for paracomplete cases with a similar proof.

\begin{thm}
Given two topological paracomplete models $M = \langle S, \sigma, V \rangle$ and $M' = \langle S', \sigma', V' \rangle$ with two continuous functions $f, f' : S \rightarrow S'$ both of which respect the valuation: $V' = f(V)=f'(V)$. If there is a homotopy $H$ between $f$ and $f$, then homotopic models satisfy the same formulae.
\end{thm}

\subsection{Classical Case}

The reason why homotopies work nicely in non-classical cases is immediate: because we stipulate that the extension of propositions to be open (or dually closed) sets. This is a strong assumption.

In the topological semantics of basic modal logic, extensions of only modal formulae are taken to be open (or closed). Can we then have results similar to those we had in non-classical case? This is the problem we are going to address in this section.

Let $N = \langle T, \eta, V \rangle$ and $N'= \langle T', \eta', V' \rangle$ be classical topological modal models. Define a homotopy $H : T \times [0,1] \rightarrow T'$. Therefore, as before, for each $t \in [0,1]$, we obtain models $N_{t} = \langle T_{t}, \eta_{t}, V_{t} \rangle$.

\begin{thm}
Given two classical topological modal models $N = \langle T, \eta, V \rangle$ and $N'= \langle T', \eta', V' \rangle$ with two continuous functions $f, f': T \rightarrow T'$ both of which respect the valuation: $V' = f(V) = f(V')$. If there is a homotopy H between $f$ and $f'$, then homotopic models satisfy the same formula.
\end{thm}

\begin{proof}
Let two classical topological modal models $N = \langle T, \eta, V \rangle$ and $N'= \langle T', \eta', V' \rangle$ with two continuous functions $g, g': T \rightarrow T'$ both of which respect the valuation: $V' = f(V) = f(V')$ be given. Let $H$ be a homotopy $H : T \times [0,1] \rightarrow T'$. We will show that homotopic models satisfy the same formula. 

Take two homotopic models $N_{t}$ and $N_{t'}$ for $t, t' \in [0, 1]$. Let $f_{t}(w_{t}) \in T_{t}$ be a point in $N_{t}$ where $f_{t} : T \rightarrow T_{t}$. Consider $N_{t}, f_{t}(w_{t}) \models \varphi$. We will show that for some $f_{t'}(w_{t'}) \in T_{t'}$, we will have $N_{t'}, f_{t'}(w_{t'}) \models \varphi$. 

Proof is by induction on the complexity of $\varphi$. First, let $\varphi = p$ for a propositional variable $p$. Then, let $N_{t}, f_{t}(w_{t}) \models p$. By using $H$, we can rewrite as $N_{t}, H(w_{t}, t) \models p$. Therefore, $H(w_{t}, t) \in V_{t}(p)$. Since $H$ is continuous on $t$, we have a continuous function from $i : t \rightarrow t'$. Now, since $V_{t} = f(V)$ we observe $H(w_{t}, t) \in f_{t}(V(p))$ which is equivalent to say $H(w_{t}, t) \in H(V(p), t)$. We can compose both sides with $i$ to get $H(w_{t'}, t') \in H(V(p), t')$ for some $w_{t'} \in T$. In short, we obtain $N_{t'}, f_{t'}(w_{t'}) \models p$.

Note that since $H$ is a homotopy, the function $i: t \rightarrow t'$ exists and is continuous. Similarly, another function $j: t' \rightarrow t$ exists and is continuous, and $j$ is needed to prove the other direction. The cases for Boolean $\varphi$ is similar and thus left to the reader.

Let us now consider the modal case $\varphi = \Box \psi$ for some $\psi$. Now, let $N_{t}, f_{t}(w_{t}) \models \Box \psi$. Thus, $f_{t}(w_{t}) \in \mathsf{Int}([\psi])$ where $[\psi]$ denotes the extension of $\psi$. Since $f_{t}$ is continuous the inverse image of an open set is open, thus $f_{t}^{-1}(\mathsf{Int}([\psi]))$ is open. For the previously constructed $i$, we observe $i^{-1} \circ (f_{t}^{-1}(\mathsf{Int}([\psi])))$ is also open as $f_{t}$ and $i$ are both continuous. Thus, the inverse image of $\mathsf{Int}([\psi])$ is open under $f_{t'}$. Therefore, by the similar reasoning, $H(w_{t'}, t') \in \mathsf{Int([\psi])}$ for some $w_{t'}$ in the neighborhood. Thus, $N_{t'}, f_{t'}(w_{t'}) \models \Box \psi$.

The reverse direction of the proof from $N_{t'}$ to $N_{t}$ is similar, and this concludes the proof.
\end{proof}

In conclusion, under a suitable valuation, isotopic models are truth invariant both in classical and non-classical cases that we have investigated.

\section{A Modal Logical Application}

Consider the following two bisimilar Kripke models $M$ and $M'$. Assume that $w, w'$ and $u, u', y'$ and $v, v', x'$ do satisfy the same propositional letters. Then it is easy to see that $w$ and $w'$ are bisimilar, and therefore satisfy the same model formulae.

\begin{figure}[htp]
\begin{picture}(200,120)
\put(40,60){\mbox{$w$}}
\put(125,12){\mbox{$v$}}
\put(125,107){\mbox{$u$}}
\put(52,62){\vector(3,2){70}}
\put(52,62){\vector(3,-2){70}}
\put(70,-5){\mbox{$M$}}
\put(200,60){\mbox{$w'$}}
\put(285,12){\mbox{$v'$}}
\put(285,107){\mbox{$u'$}}
\put(285,83){\mbox{$y'$}}
\put(285,35){\mbox{$x'$}}
\put(212,62){\vector(3,2){70}}
\put(212,62){\vector(3,-2){70}}
\put(212,62){\vector(3,1){70}}
\put(212,62){\vector(3,-1){70}}
\put(270,-5){\mbox{$M'$}}
\end{picture}
\caption{Two Bisimular Models}
\end{figure}
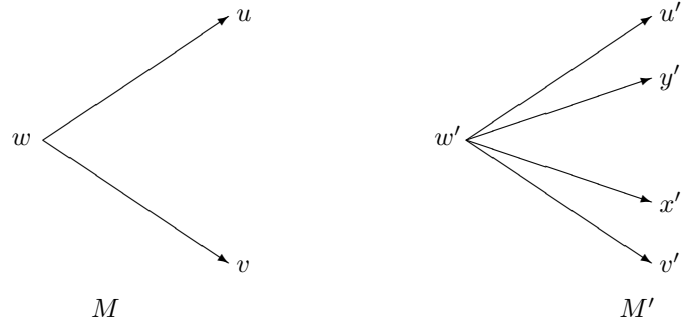

We can still pose a conceptual question about the relation between $M$ and $M'$. Even if these two models are bisimilar, they are different models. Moreover, it is plausible to \emph{contract} $M'$ to $M$ in a validity preserving fashion. Therefore, we may need to transform one model to a bisimilar model of it.
Furthermore, given a model, we may need to \emph{measure} the level of change from the fixed model to another model which is bisimilar to the given one.

Especially, in epistemic logic, such concerns do make sense. Given an epistemic situation, we can model it \emph{upto bisimulation}. In other words, from an agent's perspective, bisimilar models are indistinguishable. But, from a model theoretical perspective, they are distinguishable. Therefore, there can be several modal logical ways to model the given epistemic situation. We will now define how these models are related and different from each other by using the constructions we have presented earlier.

Before proceeding further, let us give the definition of topological bisimulations \cite{aie1}.

\begin{dfn}
Let $M= \langle S, \sigma, v \rangle$ and $M' = \langle S', \sigma', v' \rangle$ be two topological models. A topo-bisimulation is a nonempty relation $\rightleftarrows \subseteq S \times S'$ such that if $s \rightleftarrows s'$, then we have the following:
\begin{enumerate}
\item BASE CONDITION

$s \in v(p)$ if and only if $s \in v'(p)$ for any propositional variable $p$.
\item FORTH CONDITION

$s \in U \in \sigma$ implies that there exists $U' \in \sigma'$ such that $s' \in U'$ and for all $t' \in U'$ there exists $t \in U$ with $t \rightleftarrows t'$
\item BACK CONDITION
$s' \in U' \in \sigma'$ implies that there exists $U \in \sigma$ such that $s \in U$ and for all $t \in U$ there exists $t' \in U'$ with $t \rightleftarrows t'$
\end{enumerate} \end{dfn}

We can take one step further and define a homoemorphism that respect bisimulations. Let $M= \langle S, \sigma, v \rangle$ and $M' = \langle S', \sigma', v' \rangle$ be two topo-bisimular models. If there is a homoemorphism $f$ from $\langle S, \sigma \rangle$ to $\langle S, \sigma \rangle$ that respect the valuation, we call $M$ and $M$ homeo-topo-bisimilar models. We give the precise definition as follows.

\begin{dfn}
Let $M= \langle S, \sigma, v \rangle$ and $M' = \langle S', \sigma', v' \rangle$ be two topological models. A homeo-topo-bisimulation is a nonempty relation $\rightleftarrows_{f} \subseteq S \times S'$ based on a homeomorphism $f$ from $S$ into $S'$ such that if $s \rightleftarrows_{f} s'$, then we have the following:
\begin{enumerate}
\item BASE CONDITION

$s \in v(p)$ if and only if $s \in v'(p)$ for any propositional variable $p$.
\item FORTH CONDITION

$s \in U \in \sigma$ implies that there exists $f(U) \in \sigma'$ such that $s' \in f(U)$ and for all $t' \in f(U)$ there exists $t \in U$ with $t \rightleftarrows_{f} t'$
\item BACK CONDITION
$s' \in f(U) \in \sigma'$ implies that there exists $U \in \sigma$ such that $s \in U$ and for all $t \in U$ there exists $t' \in f(U)$ with $t \rightleftarrows_{f} t'$
\end{enumerate} \end{dfn}

Based on this definition, we immediately observe the following.

\begin{thm}
Homeo-topo-bisimulation preserve the validity.
\end{thm}

\begin{proof} The proof is an induction on the complexity of the formulae and thus left to the reader.\end{proof}

Notice that we can define more than one homeomorphism between topo-bisimilar models.

Now, we can discuss the homotopy of homeo-topo-bisimilar models. What we aim is the following. Given a topological model (either classical, intuitionistic or paraconsistent), we will construct two homeomorphic image of it respecting homeo-topo-bisimulation where these two homeomorphisms are homotopic. Then, by using homotopy, we will measure the level of change of the intermediate homeomorphic models with respect to these two functions.

Let $M$ be a given topological model. Construct $M_{f}$ and $M_{g}$ as the homeomorphic image of $M$ respecting the valuation where $f$ and $g$ are homeomorphism. For simplicity, assume that $M\rightleftarrows_{f} M_{f}$ and $M\rightleftarrows_{g} M_{g}$. Now, if $f$ and $g$ are homotopic, then we have functions $h_{x}$ for $x$ continuous on $[0, 1]$ with $h_{0} = f$ and $h_{1} = g$.

Therefore, given $x \in [0, 1]$ the model $M_{x}$ will be obtained by applying $h_{x}$ to $M$ respecting the valuation. Hence, $M_{0} = M_{f}$ and $M_{1} = M_{g}$. Therefore, given $M$, the distance of any homeo-topo-bisimilar model $M_{x}$ to $M$ will be $x$, and it will be the measure of non-modal change in the model. In other words, even if $M \rightleftarrows_{h(x)} M_{x}$, we will say $M$ and $M_{x}$ are $x$-different than each other.

The procedure we described offers a well-defined method of indexing the homeo-topo-bisimular models. But, indexing is not random. It is continuously on the closed unit interval.

Note that invariance results are usually used to prove undefinability results in modal logic \cite{bla1}. For example, in order to show irreflexivity is not modally definable, one needs to come up with two bisimilar models - one is irreflexive, the other is not.

Homeo-topo-bisimulations can also be used to show some topological properties are not modally definable. For instance, in this respect, dimensions of spaces is not modally definable in topological modal logic. Similarly, as trifoil and circle are homeomorphic, knots are also not definable.

\section{An Epistemic Logical Application}

Consider two believers Ann and Bob where Ann is an ordinary believer while Bob is a religious cleric of the religion that Ann is following. Therefore, they believe in the same religion and the same rules of the religion. For the sake of our example, let us assume that the religion in question is really simple and there is nothing that Ann does not know about it. In other words, even if Bob is a clergyman, Ann believes in the religion as much as Bob does. However, we feel that Bob believes it \emph{more} than Ann even though they believe in exactly the same propositions. In other words, there is still a difference between their belief. Then, what is this difference? The immediate reason for this is the fact that the \emph{extend} of Bob's knowledge is wider than that of Ann's. Namely, Bob believes \emph{better}. What does this mean?

In this context, we can ask the following two questions.
\begin{enumerate}
\item How wider is Bob's belief?
\item How is Ann's belief transformed to Bob's?
\end{enumerate}

These two questions are meaningful. Even if their language cannot tell us which one has wider knowledge, ontologically, we know that Bob has \emph{more} knowledge in some sense even if they agree on every proposition. Clearly, the reason for that is the fact that Bob considers more possible worlds for a given proposition which makes his belief more robust than Ann's.

\begin{dfn}
Given two agents $i$ and $j$ with their respective bisimilar models $M$ and $N$. We say $i$'s knowledge more robust than $j$'s if $[\varphi]^N \subseteq [\varphi]^M$ for all formula $\varphi$ in the language.
\end{dfn}

What matters is knowledge, and in basic modal logic with topological semantics, we know that the extension of modal formulas is open (or, closed dually). Therefore, in bisimilar models for Ann and Bob, the extension of Bob's robust belief is a superset of that of Ann's. Homeomorphisms and homotopies can explain this transformation from Ann's beliefs to Bob's belief.

The timestamp in the definition of the homotopy can easily be considered as a temporal parameter. In our simple example, this reading of the homotopy parameter is especially helpful. It help us to give a step by step account of the transformation between Ann's and Bob's belief.

Therefore, focusing on the transformations reveal more information on the ontology of the agents when the epistemics of the agents is already known. Robust knowledge/belief captures this notion by focusing on the ontology of the model and connects the epistemics with ontology.

\paragraph{Acknowledgement} I am grateful to Chris Mortensen, Graham Priest, Melvin Fitting  for their encouragament and comments.

\bibliographystyle{authordate3}  
\bibliography{/Users/can/Documents/Akademik/papers.bib} 

\begin{thebibliography}{}

\bibitem[\protect\citename{Aiello \& van Benthem, }2002]{aie1}
{\sc Aiello, Marco, \& van Benthem, Johan}. 2002.
\newblock A Modal Walk Through Space.
\newblock {\em Journal of Applied Non-Classical Logics}, {\bf 12}(3-4),
  319--363.

\bibitem[\protect\citename{Artemov {\em et~al.}, }1997]{art1}
{\sc Artemov, S., Davoren, J.~M., \& Nerode, A.} 1997 (June).
\newblock {\em Modal Logics and Topological Semantics for Hybrid Systems}.
\newblock Tech. rept. Mathematical Sciences Institute, Cornell University.

\bibitem[\protect\citename{Ba{\c{s}}kent, }2011]{bas3}
{\sc Ba{\c{s}}kent, Can}. 2011.
\newblock Paraconsistency and Topological Semantics.
\newblock {\em under submission, available at canbaskent.net}.

\bibitem[\protect\citename{B{\'e}ziau, }2005]{bez}
{\sc B{\'e}ziau, Jean-Yves}. 2005.
\newblock Paraconsistent Logic from a Modal Viewpoint.
\newblock {\em Journal of Applied Logic}, {\bf 3}(1), 7--14.

\bibitem[\protect\citename{Blackburn, }2000]{bla}
{\sc Blackburn, Patrick}. 2000.
\newblock Representation, reasoning and relational structures: a hybrid logic
  manifesto.
\newblock {\em Logic Journal of the IGPL}, {\bf 8}(3), 339--365.

\bibitem[\protect\citename{Blackburn {\em et~al.}, }2001]{bla1}
{\sc Blackburn, Patrick, Rijke, Maartijn~de, \& Venema, Yde}. 2001.
\newblock {\em Modal Logic}.
\newblock Cambridge Tracts in Theoretical Computer Science.
\newblock Cambridge University Press.

\bibitem[\protect\citename{Fine, }1972]{fin0}
{\sc Fine, Kit}. 1972.
\newblock In so many possible worlds.
\newblock {\em Notre Dame Journal of Formal Logic}, {\bf 13}(4), 516--520.

\bibitem[\protect\citename{Goldblatt, }2006]{gold2}
{\sc Goldblatt, Robert}. 2006.
\newblock Mathematical Modal Logic: A View of Its Evolution.
\newblock {\em In:} {\sc Gabbay, Dov~M., \& Woods, John} (eds), {\em Handbook
  of History of Logic},  vol. 6.
\newblock Elsevier.

\bibitem[\protect\citename{Goodman, }1981]{good}
{\sc Goodman, Nicolas~D.} 1981.
\newblock The Logic of Contradiction.
\newblock {\em Zeitschrift f{\"{u}}r Mathematische Logik und Grundlagen der
  Mathematik}, {\bf 27}(8-10), 119--126.

\bibitem[\protect\citename{Kremer \& Mints, }2005]{kre}
{\sc Kremer, Philip, \& Mints, Grigori}. 2005.
\newblock Dynamic Topological Logic.
\newblock {\em Annals of Pure and Applied Logic}, {\bf 131}(1-3), 133--58.

\bibitem[\protect\citename{McKinsey \& Tarski, }1944]{mck1}
{\sc McKinsey, J. C.~C., \& Tarski, Alfred}. 1944.
\newblock The Algebra of Topology.
\newblock {\em The Annals of Mathematics}, {\bf 45}(1), 141--191.

\bibitem[\protect\citename{McKinsey \& Tarski, }1946]{mck}
{\sc McKinsey, J. C.~C., \& Tarski, Alfred}. 1946.
\newblock On Closed Elements in Closure Algebras.
\newblock {\em The Annals of Mathematics}, {\bf 47}(1), 122--162.

\bibitem[\protect\citename{Mints, }2000]{min}
{\sc Mints, Grigori}. 2000.
\newblock {\em A Short Introduction to Intuitionistic Logic}.
\newblock Kluwer.

\bibitem[\protect\citename{Mortensen, }2000]{mor}
{\sc Mortensen, Chris}. 2000.
\newblock Topological Seperation Principles and Logical Theories.
\newblock {\em Synthese}, {\bf 125}(1-2), 169--178.

\bibitem[\protect\citename{Priest, }1998]{pri1}
{\sc Priest, Graham}. 1998.
\newblock What Is So Bad About Contradictions?
\newblock {\em Journal of Philosophy}, {\bf 95}(8), 410--426.

\bibitem[\protect\citename{Priest, }2002]{pri0}
{\sc Priest, Graham}. 2002.
\newblock Paraconsistent Logic.
\newblock {\em Pages  287--393 of:} {\sc Gabbay, Dov, \& Guenthner, F.} (eds),
  {\em Handbook of Philosophical Logic},  vol. 6.
\newblock Kluwer.

\bibitem[\protect\citename{Priest, }2009]{pri2}
{\sc Priest, Graham}. 2009.
\newblock Dualising Intuitionistic Negation.
\newblock {\em Principia}, {\bf 13}(3), 165--84.

\end{thebibliography}
\end{document}